\documentclass[11pt]{article}
\usepackage{latexsym}
\usepackage{times}
\usepackage{graphicx,color}
\usepackage{amsmath,amssymb, amsfonts}
\catcode`@=11
\renewcommand{\@begintheorem}[2]{\it \trivlist            
      \item[\hskip \labelsep{\bf #1\ #2{\rm :}}]}         
\renewcommand{\@opargbegintheorem}[3]{\it \trivlist       
      \item[\hskip \labelsep{\bf #1\ #2\ {\rm (#3)\/:}}]}
\def\@sect#1#2#3#4#5#6[#7]#8{\ifnum #2>\c@secnumdepth
     \def\@svsec{}\else 
     \refstepcounter{#1}\edef\@svsec{\csname the#1\endcsname{.}\hskip 1em }\fi
     \@tempskipa #5\relax
      \ifdim \@tempskipa>\z@ 
        \begingroup #6\relax
          \@hangfrom{\hskip #3\relax\@svsec}{\interlinepenalty \@M #8\par}
        \endgroup
       \csname #1mark\endcsname{#7}\addcontentsline
         {toc}{#1}{\ifnum #2>\c@secnumdepth \else
                      \protect\numberline{\csname the#1\endcsname}\fi
                    #7}\else
        \def\@svsechd{#6\hskip #3\@svsec #8\csname #1mark\endcsname
                      {#7}\addcontentsline
                           {toc}{#1}{\ifnum #2>\c@secnumdepth \else
                             \protect\numberline{\csname the#1\endcsname}\fi
                       #7}}\fi
     \@xsect{#5}}

\@addtoreset{equation}{section}
\newcommand{\Delete}[1]{}
\newcommand{\pend}{\hspace*{\fill} $\Box$}

\newtheorem{lemma}{Lemma}[section]
\newtheorem{theorem}[lemma]{Theorem}
\newtheorem{corollary}[lemma]{Corollary}
\newtheorem{proposition}[lemma]{Proposition}

\newtheorem{definition}[lemma]{Definition}

\newtheorem{remark}{Remark}

\begin{document}

\title{A Note on Ordinally Concave Functions}

\author{ Satoru FUJISHIGE,\footnote{
Research Institute for Mathematical Sciences,
Kyoto University, Kyoto 606-8502, Japan.
Email: fujishig@kurims.kyoto-u.ac.jp}\,\,\,\,
Fuhito KOJIMA,\footnote{The Department of Economics, the University of  
Tokyo, Tokyo, Japan.\ \   Email: fuhitokojima1979@gmail.com}\,\, \, and\,\,
Koji YOKOTE\,\footnote{Graduate School of Economics, the University of 
Tokyo, Japan.\,\, Email: koji.yokote@gmail.com}
}

\date{November 11, 2024}

\maketitle

\begin{abstract}

The notion of ordinal concavity of utility functions has recently been 
considered by Hafalir, Kojima, Yenmez, and Yokote in economics 
while there exist earlier related works in discrete optimization and 
operations research. 
In the present note we consider functions satisfying ordinal 
concavity and introduce a weaker notion of ordinal weak-concavity 
as well. We also investigate useful behaviors of ordinally (weak-)concave 
functions and related choice correspondences, show a characterization of 
ordinally weak-concave functions, and give an efficient algorithm for 
maximizing ordinally concave functions.
We further examine a duality in ordinally (weak-)concave functions and 
introduce the lexicographic composition of ordinally weak-concave functions.

\end{abstract}


\section{Introduction}\label{sec:1}

I.~E.~Hafalir, F.~Kojima, M.~B.~Yenmez, and K.~Yokote~\cite{HKYY2022} 
have recently considered a notion of {\it ordinal concavity} for 
utility functions defined on the set $\mathbb{Z}_{\ge 0}^E$ of nonnegative 
integer vectors with a finite nonempty set $E$. 
They have shown economic implications of ordinal 
concavity such as the path-independence property of 
choice functions associated with ordinally concave utility functions
and the rationalizability of path-independent choice rules by 
ordinally concave utility functions~\cite{YHKY2023}. 
The notion of ordinal concavity is equivalent to the one that 
was introduced by the name of 
{\it semi-strict quasi M${}^\natural$-concavity} in \cite{FarooqShioura2005}
(also see \cite{MurotaShioura2003,MurotaShioura2024} and  
 \cite{ChenLi2021}),
where the domain of the functions is the integer lattice $\mathbb{Z}^E$.

In the present note we consider {}functions satisfying ordinal 
concavity and introduce a weaker notion of ordinal weak-concavity as well.
We almost follow the notation in \cite{HKYY2022,YHKY2023}.  
Note in particular that $\emptyset$ denotes the empty set as usual 
while it also means a symbol 
that does not belong to the underlying set $E$. For any $X\in 2^E$
let $X+x=X\cup\{x\}$ for $x\in E\setminus X$ and $X-x=X\setminus\{x\}$ 
for $x\in X$. Also for $x=\emptyset$\, let $X\pm x=X$.

Let $u: 2^E\to\mathbb{R}$ be a function on the set of all subsets of a finite 
nonempty set $E$. 
The notion of {\it ordinal concavity} is defined as follows 
(see \cite{HKYY2022,YHKY2023} and 
\cite{FarooqShioura2005,MurotaShioura2003,MurotaShioura2024}).

\begin{definition}[{\bf Ordinal Concavity}] \label{def:1}
A function $u: 2^E\to \mathbb{R}$ satisfies  
{\rm ordinal concavity}\, if\, for every $X, X'\in 2^E$\, 
the following statement holds\,{\rm :} \smallskip\\
For every $x\in X\setminus X'$\, there exists 
$x'\in (X'\setminus X)\cup\{\emptyset\}$
such that 
\begin{itemize}
\item[{\rm (i)}] $u(X)<u(X-x+x')$,\, or 
\item[{\rm (ii)}] $u(X')<u(X'-x'+x)$,\, or
\item[{\rm (iii)}] $u(X)=u(X-x+x')$\, and\, $u(X')=u(X'-x'+x)$.\vspace{-0.1cm}
\end{itemize} 
\end{definition}

Let us also consider a weaker version 
which we call {\it ordinal weak-concavity} (or {\it ordinal w-concavity} 
for short) defined as follows. (See Appendix~\ref{Appendix:example3} for an example 
of a function that is ordinally weak-concave but is not ordinally 
concave.)

\begin{definition}[{\bf Ordinal Weak-Concavity}]\label{def:2} 
A function $u: 2^E\to \mathbb{R}$ satisfies  
{\rm ordinal weak-concavity} if for every $X, X'\in 2^E$ with $X\neq X'$
the following statement holds\,{\rm :} 
There exist distinct 
$x\in (X\setminus X')\cup\{\emptyset\}$ and 
$x'\in (X'\setminus X)\cup\{\emptyset\}$
such that 
\begin{itemize}
\item[{\rm (i)}] $u(X)<u(X-x+x')$,\, or
\item[{\rm (ii)}] $u(X')<u(X'-x'+x)$,\, or
\item[{\rm (iii)}] $u(X)=u(X-x+x')$\, and\, $u(X')=u(X'-x'+x)$.
\end{itemize}
\end{definition}

\begin{remark}{\rm 
A notion of {\it weak semi-strict quasi M-concavity} (denoted by
(SSQM${}_w$)) is considered in \cite{MurotaShioura2003}. 
Just as we can obtain 
an M${}^\natural$-convex function from an M-convex function by a projection 
of the domain $\mathbb{R}^n$ along an axis into a one-dimension-lower 
coordinate space $\mathbb{R}^{n-1}$
(see \cite{Fuji2005,Murota2003,MurotaShioura1999,ShiouraTamura2015}), 
we can define a notion of 
{\it weak semi-strict quasi M${}^\natural$-concavity} from weak 
semi-strict quasi M-concavity, which has not been explicitly considered  
in the literature. The notion of ordinal weak-concavity
given above is a set-theoretical version of \lq weak semi-strict 
quasi M${}^\natural$-concavity.'
\pend}
\end{remark}

The present note is organized as follows. 
 In Section~\ref{sec:wconcavity} we examine some useful behaviors of 
ordinally w-concave functions and 
show a characterization of ordinally w-concave functions.
We also propose an efficient algorithm for maximizing ordinally concave 
functions. 
In Sections~\ref{sec:choice2} and \ref{sec:UM} we discuss 
path-independent choice functions
for ordinally concave functions \cite[Theorem 2]{HKYY2022} in view of 
choice correspondences and we also 
examine behaviors of choice functions associated with ordinally 
w-concave functions and related choice correspondences 
in Section~\ref{sec:Umax}.
In Section~\ref{sec:discussions} we examine a duality property in 
ordinal concavity and 
introduce the lexicographic composition of two ordinally w-concave 
functions. Section~\ref{sec:concluding} gives concluding remarks.

\section{Ordinally Weak-Concave Functions}\label{sec:wconcavity}

For any $X, Y\in 2^E$ define $X\Delta Y=(X\setminus Y)\cup(Y\setminus X)$ 
(the symmetric difference of $X$ and $Y$). 
Also for any finite set $X$ denote by $|X|$ the number of elements of $X$.

\subsection{Fundamental operations on functions on $2^E$}
\label{sec:w-concavity0}

Consider any function $u: 2^E\to\mathbb{R}$. 
For any nonempty $X\subseteq E$ define $u^X: 2^X\to\mathbb{R}$ by
\begin{equation}\label{eq:fop1}
  u^X(Z)=u(Z) \qquad (\forall Z\in 2^X).
\end{equation}
We call $u^X$ the {\it reduction} of $u$ by $X$ 
(or the {\it restriction} of $u$ on $X$). 
Also for any $X\subset E$ define  $u_X: 2^{E\setminus X}\to\mathbb{R}$ by
\begin{equation}\label{eq:fop2}
  u_X(Z)=u(Z\cup X)-u(X) \qquad (\forall Z\in 2^{E\setminus X}).
\end{equation}
We call $u_X$ the {\it contraction} of $u$ by $X$.
Moreover, for any $X,Y\in 2^E$ with $X\subset Y$ define 
$u_X^Y: 2^{Y\setminus X}\to\mathbb{R}$ by
\begin{equation}\label{eq:fop3}
  u_X^Y(Z)=u(Z\cup X)-u(X) \qquad (\forall Z\in 2^{Y\setminus X}).
\end{equation}
We call $u_X^Y$ a {\it minor} of $u$ obtained by the reduction by $Y$ and 
then by the contraction by $X$.
It should be noted that 
\begin{itemize}
\item the operations of reduction and contraction keep ordinal 
(w-)concavity and hence every minor of an ordinally (w-)concave $u$ is 
ordinally (w-)concave.
\end{itemize}

\subsection{Fundamental properties of ordinally weak-concave functions}
\label{sec:w-concavity}

Let us consider the following simultaneous exchange property $(\dagger)$
for a nonempty family $\mathcal{F}$ of subsets of $E$.
\begin{itemize}
\item[$(\dagger)$] For every $X, X'\in \mathcal{F}$ with $X\neq X'$
there exist distinct 
$x\in (X\setminus X')\cup\{\emptyset\}$ and 
$x'\in (X'\setminus X)\cup\{\emptyset\}$ such that 
$X-x+x',\, X'+x-x'\in\mathcal{F}$.
\end{itemize}

We have the following proposition. 
Although this proposition appears to be implicitly known in the literature,
we give a proof of it for completeness. For any $X\subseteq E$ denote by
$\chi_X$ the characteristic vector of $X$, i.e., $\chi_X\in \mathbb{R}^E$
and $\chi_X(i)=1$ for all $i\in X$ and $\chi_X(i)=0$ for all 
$i\in E\setminus X$.
Also we write $\chi_{\{x\}}$ simply as $\chi_x$ for any $x\in E$.

\begin{proposition}\label{prop:exc0}
Let $\mathcal{F}$ be any nonempty family of subsets of $E$. 
Then $\mathcal{F}$ satisfies the simultaneous exchange property 
$(\dagger)$ if and only if $\mathcal{F}$ forms an 
M${}^\natural$-convex set {\rm (}or a generalized matroid\,{\rm )} 
on $E$. 
\medskip\\
{\rm (Proof) {\sl The if part}\/: Suppose that $\mathcal{F}$ forms an 
M${}^\natural$-convex set on $E$.
Consider any distinct $X, X'\in\mathcal{F}$. Because of the symmetry 
between $X$ and $X'$ in $(\dagger)$ we can suppose without loss of generality 
that $X\setminus X'\neq \emptyset$.
Then for any $x\in X\setminus X'$ there exists 
$x'\in (X'\setminus X)\cup\{\emptyset\}$ such that 
$X-x+x', X'+x-x'\in\mathcal{F}$, due to M${}^\natural$-convexity 
of $\mathcal{F}$ (see \cite{Murota2003}). 
Hence the if part of $(\dagger)$ holds.

{\sl The only-if part}\/: 
For a nonempty family $\mathcal{F}$ of subsets of $E$ that satisfies 
the simultaneous exchange property $(\dagger)$, 
consider the convex hull, denoted by ${\rm Conv}(\mathcal{F})$, of the set of 
characteristic vectors $\chi_X$ of all $X\in\mathcal{F}$. 
Note that $\chi_X$ for every $X\in\mathcal{F}$ is an extreme point of 
the polytope ${\rm Conv}(\mathcal{F})$ and we have 
\begin{equation}\label{eq:convZ}
  {\rm Conv}(\mathcal{F})\cap \mathbb{Z}^E =\{\chi_X\mid X\in\mathcal{F}\}.
\end{equation}
Moreover, it follows from $(\dagger)$ that for any pair
of adjacent extreme points $\chi_X$ and $\chi_{X'}$ of 
${\rm Conv}(\mathcal{F})$ there exist distinct 
$x\in (X\setminus X')\cup\{\emptyset\}$ and 
$x'\in (X'\setminus X)\cup\{\emptyset\}$ such that 
$X-x+x',\, X'+x-x'\in\mathcal{F}$. 
Since the direction vectors from $\chi_X$ to $\chi_{X-x+x'}$ and 
from $\chi_{X'}$ to $\chi_{X'+x-x'}$ are given by
\begin{equation}\label{eq:prop0}
\chi_{X-x+x'}-\chi_X=\chi_{x'}-\chi_x,\quad 
\chi_{X'+x-x'}-\chi_{X'}=\chi_{x}-\chi_{x'}
\end{equation} 
and since $\chi_X$ and $\chi_{X'}$ are adjacent
extreme points of ${\rm Conv}(\mathcal{F})$, it follows that
the edge direction between $\chi_X$ and $\chi_{X'}$ is a 
non-zero scalar multiple of $\chi_{x'}-\chi_x$. (Actually we have
$\chi_{X'}-\chi_X=\chi_{x'}-\chi_x$ since there is no integral point
in the open line-segment between $\chi_X$ and $\chi_{X'}$.)
This fact together with (\ref{eq:convZ}) implies that $\mathcal{F}$ forms an 
M${}^\natural$-convex set on $E$ (due to, e.g., \cite[Theorem 17.1]{Fuji2005}).
\pend}
\end{proposition}

Now, for any function $u: 2^E\to\mathbb{R}$ define 
\begin{equation}\label{eq:f1}
 {\bf D}_u^*={\rm Arg}\max\{u(X)\mid X\subseteq E\}.
\end{equation}
We see the following fact.

\begin{lemma}\label{lem:f0}
For every ordinally w-concave function $u$   
the family ${\bf D}_u^*$ of subsets of $E$ satisfies the property $(\dagger)$ 
with $\mathcal{F}={\bf D}_u^*$.
\medskip\\
{\rm (Proof) For every ordinally w-concave function $u$ 
we see that for every $X, X'\in {\bf D}_u^*$ 
with $X\neq X'$ there exist distinct 
$x\in (X\setminus X')\cup\{\emptyset\}$ and 
$x'\in (X'\setminus X)\cup\{\emptyset\}$ such that only Condition (iii) of
Definition~\ref{def:2} holds because of the definition of ${\bf D}_u^*$. 
Hence $\mathcal{F}={\bf D}_u^*$ satisfies the property $(\dagger)$.
\pend}
\end{lemma}

From Lemma~\ref{lem:f0} and the only-if part of Proposition~\ref{prop:exc0} 
we have the following lemma.

\begin{lemma}\label{lem:f1}
For every ordinally w-concave function $u$   
the set ${\bf D}_u^*$ given by {\rm (\ref{eq:f1})}\, forms an 
M${}^\natural$-convex set 
on $E$.
\end{lemma}

Next, we show other fundamental facts about ordinally w-concave functions.
For every $X\in 2^ E$ denote by $X^*$ a set $Z\in {\bf D}_u^*$ 
 that minimizes $|X\Delta Z|$ among ${\bf D}_u^*$. 
(In this note $X'$ is used as a variable set that is independent of $X$, 
while for a given $X\in 2^E$, $X^*$ is a member of 
${\rm Arg}\min\{|X\Delta Z|\mid Z\in {\bf D}_u^*\}$.)

\begin{theorem}\label{th:main}
Let $u$ be any ordinally w-concave function. 
Consider any $X, Y\in 2^E\setminus {\bf D}_u^*$  satisfying 
\begin{equation}\label{eq:00}
 X\cap X^*\subseteq Y\subseteq X\cup X^*.
\end{equation}
Then there exist distinct 
$x\in (X^*\setminus Y)\cup\{\emptyset\}$ and 
$y\in (Y\setminus X^*)\cup\{\emptyset\}$
such that  $u(Y)<u(Y-y+x)$.
\medskip\\
{\rm (Proof)
Choose any $X, Y\in 2^E\setminus {\bf D}_u^*$ that satisfy (\ref{eq:00}).
Note that $Y\neq X^*$ since $Y\notin {\bf D}_u^*$. 
It follows from ordinal w-concavity of $u$ and the definition of $X^*$
that for $X\gets Y$ and $X'\gets X^*$ only (i) in the definition of
ordinal w-concavity holds.
\pend}
\end{theorem}

\begin{corollary}\label{cor:main}
Let $u: 2^E\to\mathbb{R}$ be any ordinally w-concave function.
Consider any $Y\in 2^E\setminus {\bf D}_u^*$ and $Z\in {\bf D}_u^*$.
Then there exist distinct $x\in (Z\setminus Y)\cup\{\emptyset\}$ and 
$y\in(Y\setminus Z)\cup\{\emptyset\}$ such that $u(Y)<u(Y-y+x)$.
\medskip\\
{\rm (Proof) For any $Y\in 2^E\setminus {\bf D}_u^*$ and 
$Z\in {\bf D}_u^*$ consider 
the minor $\bar{u}\equiv u_{Y\cap Z}^{Y\cup Z}$ of $u$, which is ordinally 
w-concave.
Note that $(Y\cup Z)\setminus (Y\cap Z)\neq\emptyset$. 
Let us use the unary operator $(\cdot)^*$
(defined for $u$) for $\bar{u}$ as well. 
Note that since $Y\notin {\bf D}_u^*$ and  $Z\in {\bf D}_u^*$, we have 
$(Y\setminus Z)^*\neq Y\setminus Z$ for $\bar{u}$.
Then, considering $u\gets\bar{u}$, 
$X\gets Y\setminus Z$ and $X^*\gets (Y\setminus Z)^*$ (for 
$\bar{u}$) in Theorem~\ref{th:main} with $X=Y$, we have
\begin{itemize}
\item there exist distinct $x\in((Y\setminus Z)^*\setminus Y)\cup\{\emptyset\}$
and $y\in (Y\setminus (Y\setminus Z)^*)\cup\{\emptyset\}$ such that 
$\bar{u}(Y\setminus(Y\cap Z))<\bar{u}((Y-y+x)\setminus (Y\cap Z))$.
\end{itemize}
Since $(Y\setminus Z)^*\setminus Y\subseteq Z\setminus Y$ 
 and $Y\setminus (Y\setminus Z)^*\subseteq Y\setminus Z$
and since 
$\bar{u}((Y-y+x)\setminus (Y\cap Z))=u(Y-y+x)-u(Y\cap Z)$ and
$\bar{u}(Y\setminus(Y\cap Z))=u(Y)-u(Y\cap Z)$, 
the pair of $x$ and $y$ is a desired one 
for the present corollary.
\pend}
\end{corollary}

{From} Theorem~\ref{th:main} we see the following corollaries,
where we suppose that $u$ satisfies ordinal w-concavity.

\begin{corollary}\label{cor:f1}
For any $X\in 2^E\setminus {\bf D}_u^*$ there exists a sequence of distinct 
subsets $Y_0(=X), Y_1, \cdots, Y_k$ for a positive integer $k\le |E|$ 
such that 
\begin{itemize}
\item[{\rm 1.}]\, $Y_i\in 2^E\setminus {\bf D}_u^*$\, 
for each $i\in\{0,1,\cdots,k-1\}$\,,
\item[{\rm 2.}]\, $u(Y_0)<u(Y_1)<\cdots <u(Y_k)$ with $Y_k\in {\bf D}_u^*$\,,
\item[{\rm 3.}]\, for each $i\in\{0,1,\cdots,k-1\}$ we have 
$Y_{i+1}=Y_i-y_i+x_i$ for  
distinct  $y_i\in (Y_i\setminus X^*)\cup\{\emptyset\}$ and 
  $x_i\in (X^*\setminus Y_i)\cup\{\emptyset\}$.
\end{itemize}
{\rm (Proof) It follows from Theorem~\ref{th:main} that repeating the 
transformation $Y\gets Y-y+x$ as far as $Y\notin {\bf D}_u^*$, we obtain   
 $Y\in {\bf D}_u^*$ after at most $|X\Delta X^*|$ such transformations, 
each increasing the value of $u(Y)$.
\pend}
\end{corollary}

A simple consequence of Theorem~\ref{th:main} is also given as follows. 
For any linear ordering $L=(e_1,e_2,\cdots,e_k)$ 
of distinct $k$ elements of $E$ with $k\le |E|$ define $L_i=\{e_1,\cdots,e_i\}$
(the set of the initial $i$ elements of $L$) for each $i\in\{1,\cdots,k\}$.

\begin{corollary}\label{cor:f0} 
For any $Z\in {\bf D}_u^*$ of minimum cardinality $|Z|$ 
there exists a linear ordering $L=(e_1,\cdots,e_k)$ 
of elements of $Z$ such that 
$u(\emptyset)<u(L_1)<\cdots<u(L_k)$ with $L_k=Z$. 
\medskip\\
{\rm (Proof)
Starting from $X=\emptyset$, we can reach any maximizer 
$Z$ of $u$ having the minimum cardinality by
repeating the transformation of $Y$ only in a form of $Y+x$ for some 
$x\in Z\setminus Y$ as in Corollary~\ref{cor:f1}. 
Consider $X\gets Y$ and $X'\gets Z$ in the definition of ordinal w-concavity
and note that $X\subset X'$ and $X'\in {\bf D}_u^*$ having the minimum 
cardinality.
\pend}
\end{corollary}

For any $X\in 2^E$ define the {\it neighborhood} ${\bf N}(X)$ of $X$ by
\begin{equation}\label{eq:N}
{\bf N}(X)=\{X-x+x'\mid x\in X\cup\{\emptyset\}, 
x'\in (E\setminus X)\cup\{\emptyset\}\}.
\end{equation}
Note that $|{\bf N}(X)|$ is ${\rm O}(|E|^2)$ for any $X\in 2^E$.
The following corollary is a projected version of 
\cite[Theorem 4.2(ii)]{MurotaShioura2003}.

\begin{corollary}\label{cor:f2}
A set $X\in 2^ E$ is a maximizer of $u$ if and only if $X$ attains 
the maximum of $u(Z)$ among all $Z\in {\bf N}(X)$. 
\medskip\\
{\rm (Proof) It suffices to show the if part.
Suppose that $X$ is not a maximizer of $u$. 
Then it follows from Corollary~\ref{cor:main} that 
there exists $Z\in {\bf N}(X)$ such that $u(Z)>u(X)$.
This completes the proof of the present corollary.
\pend}
\end{corollary}

Corollary~\ref{cor:f2} leads us to a simple hill-climbing algorithm 
to maximize $u$ satisfying ordinal w-concavity as follows 
(cf.~\cite[Sec.~4.2]{MurotaShioura2003}).  
\medskip\\
{\sf Algorithm 1} \smallskip\\
 {\bf Step 0}:  
Choose any $X\in 2^E$ and put $Y\gets X$ ; \\
{\bf Step 1}: 
While there exists $Z\in {\bf N}(Y)$ such that $u(Z)>u(Y)$,\, 
do the following:  \\
\qquad \qquad \quad Choose any $Z\in {\bf N}(Y)$ such that $u(Z)>u(Y)$\,;\\
\qquad \qquad \quad Put $Y\gets Z$\,; \\
{\bf Step 2}: Return $Y$ ; \medskip

Trivially, {\sf Algorithm 1} terminates after updating $Y$ in {\bf Step 1}
at most $\nu$ times, where $\nu$ is the number of 
distinct function values $u(X)$ for all $X\in 2^E$. 

Hafalir et al.~\cite[Theorem 1]{HKYY2022} show that when $u$ satisfies 
ordinal concavity, we can find a maximizer of $u$ after ${\rm O}(|E|^3)$ 
updates of $Y$ in {\bf Step 1} by appropriately choosing $Z$ 
(also see \cite{MurotaShioura2003,MurotaShioura2024}).
It is an interesting problem 
to find an algorithm faster than {\sf Algorithm~1} 
(if any) for ordinally w-concave functions. 
In Section~\ref{sec:max} we show an algorithm faster than the one given by 
Hafalir et al.~\cite{HKYY2022} for ordinally concave functions.

\subsection{A characterization of ordinal weak-concavity}
\label{sec:characterization}

For any $X, Y\in 2^E$ with $X\subseteq Y$ define 
$[X,Y]\equiv\{Z\in 2^E\mid X\subseteq Z\subseteq Y\}$ 
(an {\it interval} in $2^E$) and 
\begin{equation}\label{eq:ff1}
 {\bf C}_u(X,Y)={\rm Arg}\max\{u(Z)\mid Z\in [X,Y]\}.
\end{equation}

\begin{lemma}\label{lem:ff1}
Let $u: 2^E\to\mathbb{R}$ be any function satisfying ordinal 
w-concavity.
Then for every $X, Y\in 2^E$ and every $Z\in[X\cap Y,X\cup Y]$,\, if\,  
$Z$ maximizes $u$ over 
\[{\bf N}(Z)\cap[X\cap Y,X\cup Y],\] 
then $Z$ maximizes
 $u$ over \[[X\cap Y,X\cup Y],\] 
{\rm i.e.}, $Z\in{\bf C}_u(X\cap Y,X\cup Y)$.
\medskip\\
{\rm (Proof) Suppose that $u: 2^E\to\mathbb{R}$ is a function 
satisfying ordinal w-concavity. 
Consider any $X, Y\in 2^E$. 
We can suppose that $X\neq Y$. Then the present lemma is equivalent to 
the statement of Corollary~\ref{cor:f2} using the minor $u^{X\cup Y}_{X\cap Y}$
in place of $u$. Recall that $u^{X\cup Y}_{X\cap Y}$ is ordinally w-concave.
\pend}
\end{lemma}

Also, similarly as Lemma~\ref{lem:f1} we can show the following lemma.

\begin{lemma}\label{lem:XY}
Let $u: 2^E\to\mathbb{R}$ be any function satisfying ordinal 
w-concavity.
For any $X, Y\in 2^E$ with $X\subset Y$ 
the set ${\bf C}_u(X,Y)$ forms an M${}^\natural$-convex set.
\medskip\\
{\rm (Proof) For any $X, Y\in 2^E$ with $X\subset Y$
consider the minor $u^Y_X$ of $u$ and apply Lemma~\ref{lem:f1} 
for $u^Y_X$ in place of $u$. 
\pend}
\end{lemma}

Now, we show the following theorem, which means that the properties of 
$u$ shown in Lemmas~\ref{lem:ff1} and \ref{lem:XY} actually characterize 
ordinal w-concavity.

\begin{theorem}\label{th:ff1}
A function $u: 2^E\to\mathbb{R}$ satisfies ordinal w-concavity 
if and only if the following two statements hold\,{\rm :}
\begin{itemize}
\item[{\bf (M)}]\, For any $X, Y\in 2^E$ with $X\subset Y$ 
the set ${\bf C}_u(X,Y)$ forms an M${}^\natural$-convex set.
\item[{\bf (N)}]\, 
For every $X, Y\in 2^E$ and every $Z\in[X\cap Y,X\cup Y]$, if  
$Z$ maximizes $u$ over ${\bf N}(Z)\cap[X\cap Y,X\cup Y]$, 
then $Z$ maximizes $u$ over $[X\cap Y,X\cup Y]$.
\end{itemize}
{\rm (Proof)  The only-if part follows from Lemmas~\ref{lem:ff1} and 
\ref{lem:XY}. 
We show the if part.

Suppose that {\bf (M)} and {\bf (N)} hold. Consider any $X, Y\in 2^E$ 
with $X\neq Y$. 
Then we have the following three cases:
\begin{enumerate}
\item $X\notin{\bf C}_u(X\cap Y,X\cup Y)$,\, or
\item $Y\notin{\bf C}_u(X\cap Y,X\cup Y)$,\, or
\item $X\in{\bf C}_u(X\cap Y,X\cup Y)$\, and\, 
$Y\in{\bf C}_u(X\cap Y,X\cup Y)$.
\end{enumerate}
It follows from {\bf (N)} that in Case 1 or Case 2,  
for some distinct $x\in (X\setminus Y)\cup\{\emptyset\}$ and 
$y\in (Y\setminus X)\cup\{\emptyset\}$ 
we have
\begin{itemize}
\item[(i)]  $u(X-x+y)>u(X)$,\, or
\item[(ii)] $u(Y+x-y)>u(Y)$.
\end{itemize}
Moreover, it follows from {\bf (M)} that in Case 3, for some distinct 
$x\in (X\setminus Y)\cup\{\emptyset\}$ and 
$y\in (Y\setminus X)\cup\{\emptyset\}$ we have \smallskip\\
\ \ \ (iii)\, $u(X-x+y)=u(X)$\, and\, $u(Y+x-y)=u(Y)$.

This completes the proof of the ordinal w-concavity of $u$.
\pend}
\end{theorem}

\subsection{Ordinal concavity vs.~ordinal weak-concavity}
\label{sec:strengthening}

The lemmas and theorems shown in 
Section~\ref{sec:w-concavity} can be 
strengthened if we consider functions satisfying 
ordinal concavity instead of ordinal weak-concavity.
For example, Corollary~\ref{cor:f1} is strengthened as follows.
Suppose that $u: 2^E\to\mathbb{R}$ is a function satisfying 
ordinal concavity.

\begin{corollary}\label{cor:f1a}
For any $X\in 2^E\setminus {\bf D}_u^*$ there exists a sequence of distinct 
subsets $Y_0(=X), Y_1, \cdots, Y_k$ for a positive integer $k\le |E|$ 
such that
\begin{itemize}
\item[{\rm 1.}]\, $Y_i\in 2^E\setminus {\bf D}_u^*$\,\, for each\, 
$i\in\{0,1,\cdots,k-1\}$\,,
\item[{\rm 2.}]\, $u(Y_0)<u(Y_1)<\cdots <u(Y_k)$\, with \, 
$Y_k\in {\bf D}_u^*$\,,
\item[{\rm 3.}]\, for some integer $\ell$ with $0\le \ell\le k$ we have
\begin{itemize}
\item[{\rm (a)}] for each $i\in\{0,1,\cdots,\ell-1\}$, 
$Y_{i+1}=Y_i-x_i+x_i'$ such that 
  $x_i\in Y_i\setminus Y_{i+1}$ and 
  $x_i'\in (Y_{i+1}\setminus Y_i)\cup\{\emptyset\}$ and
\item[{\rm (b)}] for each $i\in\{\ell,\cdots,k-1\}$, $Y_{i+1}=Y_i+x_i'$ 
such that $x_i'\in Y_{i+1}\setminus Y_i$\,.
\end{itemize}
\end{itemize}
\end{corollary}

Also, Corollary~\ref{cor:f0} is strengthened as follows.
 
\begin{corollary}\label{cor:f0a} 
Consider any $X\in {\bf D}_u^*$ of minimum cardinality. Then 
for every linear ordering $L=(e_1,\cdots,e_k)$ 
of elements of $X$ we have
$u(\emptyset)<u(L_1)<\cdots<u(L_k)$ with $L_k=X$. 
\end{corollary}

Corollaries~\ref{cor:f0}, \ref{cor:f1a}, and \ref{cor:f0a} reflect 
combinatorial 
structures like anti\-matroids or convex geometries behind functions
satisfying ordinal {(w-)concavity} (cf.~\cite{Fuji2015,Koshevoy1999,MR2001}).

\subsubsection{Maximizing ordinally concave functions}\label{sec:max}

Now let us consider an algorithm for maximizing an ordinally 
concave function $u: 2^E\to\mathbb{R}$. Recall that a set $X\in 2^E$ is 
called a local 
maximizer of $u$ if $X$ maximizes $u$ in the neighborhood ${\bf N}(X)$ 
defined by (\ref{eq:N}). 
Here note that we consider ordinally concave functions 
but not ordinally w-concave functions in general.
\medskip

\noindent
{\sf Algorithm 2}\smallskip\\
{\bf Input}: An ordinally concave function $u: 2^E\to\mathbb{R}$;\\
{\bf Step 0}: Put $W\gets\emptyset$ and let $u'=u$;\\
{\bf Step 1}: While $\emptyset$ is not a local maximizer of $u'$, 
do the following:\\
\quad Choose any  
$x^*\in {\rm Arg}\max\{u(\{x\})\mid x\in E\setminus W\}$;\,
$W\gets W\cup\{x^*\}$;\\
\quad Let $u'$ be the contraction $u_W$ of the original $u$ 
by the updated $W$;\\
{\bf Step 2}: Return $W$;

\begin{theorem}\label{th:1}
For any ordinally concave function $u: 2^E\to\mathbb{R}$ 
{\sf Algorithm 2} finds a maximizer of $u$ and requires ${\rm O}(|E|^2)$
function calls for $u$ in total.
\medskip\\
{\rm (Proof) 
If the empty set $\emptyset$ is a local maximizer of $u$,
then it is the global maximizer due to Corollary~\ref{cor:f2}. 
Hence suppose that 
$\emptyset$ is not a local maximizer of $u$ and choose any
$x^*\in {\rm Arg}\max\{u(\{x\})\mid x\in E\setminus W\}$ with 
$W=\emptyset$ initially. 

We show that there exists some global maximizer $U\in{\bf C}_u(E)$ 
such that $x^*\in U$. 
For a given $U\in{\bf C}_u(E)$ suppose that $x^*\notin U$.
Then, putting $X\gets\{x^*\}$ and $X'\gets U$, 
for these $X$ and $X'$ we have that for any $x\in X\setminus X'=\{x^*\}$ 
there exists
$x'\in(X'\setminus X)\cup\{\emptyset\}$ 
such that 
one of the conditions (i), (ii), and (iii) for ordinal concavity holds.
That is, 
\begin{itemize}
\item[{\rm (i)}] $u(X)<u(X-x+x')$,\, or 
\item[{\rm (ii)}] $u(X')<u(X'-x'+x)$,\, or
\item[{\rm (iii)}] $u(X)=u(X-x+x')$\, and\, $u(X')=u(X'-x'+x)$.\vspace{-0.1cm}
\end{itemize} 
Note that $x=x^*$ and 
$X-x^*+x^{\prime}=x^{\prime} \in U \cup\{\emptyset\}$. 
Hence condition (i) does not hold because of the definition of $x^*$ and 
since $u\left(\{x^*\}\right)>u(\emptyset)$. Condition (ii)
is null because of the definition of $U$. It follows that Condition (iii) 
holds. 
Then $X'-x'+x^*$ is a global maximizer and satisfies $x^*\in X'-x'+x^*$. 

Consequently, updating $W\gets W\cup\{x^*\}$ and considering the 
contraction $u'=u_W$ of $u$ by $W(=\{x^*\}\, \text{currently})$, 
the updated $u'$ is also ordinally concave, and moreover, 
for any maximizer $U'$ of  $u'=u_W$ we have a maximizer $U'\cup W$ of $u$.

Repeating this process until $\emptyset$ becomes a local maximizer of updated 
$u'=u_W$, the finally obtained $W$ is a global maximizer of $u$.

We see that each $x^*$ in {\bf Step 1} is computed by at most $|E|$ 
function calls for $u$ and the While loop of {\bf Step 1} is repeated 
at most $|E|$ times. Hence the algorithm terminates after ${\rm O}(|E|^2)$
function calls for $u$.
\pend}
\end{theorem}

\begin{remark}{\rm 
A crucial point in {\sf Algorithm 2} is that we start from the special
set $\emptyset$ and employ the operation of contraction to execute 
the algorithm in a sort of recursive way. It can also be understood as an 
effective use of \lq maximizer-cut property' 
for ordinally concave functions  shown
in \cite{HKYY2022,MurotaShioura2003,MurotaShioura2024} that treat more
general problem setting where the empty set $\emptyset$ may not belong 
to the effective domain of $u$ and the initial starting point can be 
arbitrarily chosen.
If we apply {\sf Algorithm 2} to an M${}^\natural$-concave function 
$g: {\bf Q}\to\mathbb{R}$ with $\emptyset\in {\bf Q}\subseteq 2^E$, 
a special  ordinally concave function, 
then it becomes an upward-steepest ascent algorithm for such a function  
(see \cite[Theorem~3.7]{ShiouraTamura2015}).
\pend}
\end{remark}

\section{Choice functions and choice correspondences}\label{sec:choice}

Let $u: 2^E\to\mathbb{R}$ be a (utility) function. 
Any function $C: 2^E\to 2^E$ satisfying $C(X)\subseteq X$ $(\forall X\in 2^E)$ 
is called a {\it choice function} on $2^E$. 
For each $X\in 2^E$ define 
${\bf C}_u(X)={\rm Arg}\max\{u(Y)\mid Y\subseteq X\}$. 
Note that ${\bf C}_u(X)={\bf C}_u(\emptyset,X)$ defined in the previous 
section.
The mapping ${\bf C}_u: 2^E\to 2^{2^E}$ is called the 
{\it choice correspondence rationalized by} $u$. 
If $C$ satisfies $C(X)\in{\bf C}_u(X)$ for each $X\in 2^E$, 
we call such $C$ a {\it choice function associated with} $u$.

\subsection{A choice function associated with an ordinally concave function}
\label{sec:choice2}



The results obtained in the previous section lead us to  
the following theorem, which was shown by Farooq and 
Shioura~\cite{FarooqShioura2005}. 

\begin{theorem}[{\cite[Theorem~4.1]{FarooqShioura2005}}]\label{th:w-choice}
Let $u: 2^E\to\mathbb{R}$ be an ordinally concave function, 
${\bf C}_u$ be the choice correspondence rationalized by $u$, and 
$C: 2^E\to 2^E$ be a choice function associated with $u$.
Then the following two statements {\rm (I)} and {\rm (II)} hold\,{\rm : }
\begin{itemize}
\item[{\rm (I)}] For every $X, Y\in 2^E$ and every $U\in{\bf C}_u(X)$ there 
exists $Z\in{\bf C}_u(X\cup Y)$ such that $Z\cap X\subseteq U$.
In other words, for every $X, Y\in 2^E$  we have 
\[{\bf C}_u(X\cup Y)\cap{\bf C}_u(C(X)\cup (Y\setminus X))\neq\emptyset.\] 
\item[{\rm (II)}]  For every $X, Y \in 2^E$ and every 
$Z \in \mathbf{C}_u(X\cup Y)$ 
there exists $U \in \mathbf{C}_u(X)$ such that 
$Z\cap X\subseteq U$.
In other words, for every $X, Y\in 2^E$ we have 
\[{\bf C}_u(X)\cap{\bf C}_u(C(X\cup Y)\cap X,X)\neq\emptyset.\] 
\end{itemize}
\end{theorem}

\begin{remark}{\rm 
The above theorem, Theorem~\ref{th:w-choice}, does not necessarily hold for 
functions satisfying ordinal weak-concavity but not ordinal 
concavity. See an example shown in Appendix~\ref{Appendix:example3}.
\pend}
\end{remark}

\begin{remark}{\rm 
Under a stronger assumption that $u$ is an M${}^\natural$-concave function, 
Murota \cite[Theorem~3.8]{Murota2016} showed (I) and (II) in the above theorem.
The properties of (I) and (II) are known as the substitutability of 
the choice correspondence ${\bf C}_u$ (see \cite{Sotomayor1999}).
The substitutability plays a crucial role in the two-sided matching 
setting (see, e.g., \cite{Fleiner2003,HatfieldMilgrom2005,Sotomayor1999}) 
in that 
a stable matching exists under substitutability and Sen's $\alpha$:
$\forall X\in 2^E, \forall Y\in {\bf C}_u(X):\, 
Y\subset X'\subset X \Longrightarrow Y\in {\bf C}_u(X')$ 
(see \cite{CheKimKojima2021}). Note that 
Sen's $\alpha$ holds in our problem setting because the choice correspondence 
is rationalized by $u$.
\pend}
\end{remark}

\subsection{The unique-maximizer condition}\label{sec:UM}

Let us assume that $u$ satisfies 
the following {\it unique-maximizer} condition {\bf (UM)} in addition 
to ordinal concavity.
\begin{itemize}
\item[{\bf  (UM)}] For every $X\subseteq E$ there uniquely exists a 
maximizer of 
 $\max\{u(Y)\mid Y\subseteq X\}$, i.e., $|{\bf C}_u(X)|=1$. 
\end{itemize}
Then the choice function $C: 2^E\to 2^E$ associated with $u$ is uniquely 
determined, i.e., ${\bf C}_u(X)=\{C(X)\}$ for all $X\in 2^E$. Hence 
we can identify ${\bf C}_u(X)$ with $C(X)$  for all $X\in 2^E$.

The following theorem was shown in \cite{HKYY2022} for functions
satisfying ordinal concavity.

\begin{theorem}[{\cite[Theorem 2]{HKYY2022}}]\label{th:HKYY2022}
For any function $u: 2^E\to \mathbb{R}$ satisfying ordinal concavity and 
the unique-maximizer 
condition {\bf (UM)}, the choice function $C: 2^E\to 2^E$ associated with $u$ 
is path-independent, i.e., for every $X, Y\in 2^E$ 
we have $C(X\cup Y)=C(C(X)\cup Y)$.
\end{theorem}

Theorem~\ref{th:w-choice} actually leads us to the following theorem 
for functions satisfying ordinal concavity and Condition {\bf (UM)}.

\begin{theorem}\label{th:choice2}
For any function $u: 2^E\to \mathbb{R}$ satisfying ordinal concavity 
and the unique-maximizer condition {\bf (UM)}, the choice function 
$C: 2^E\to 2^E$ associated with $u$ satisfies  
$C(X\cup Y)=C(C(X)\cup (Y\setminus X))$ for all $X, Y\in 2^E$. 
\medskip\\
{\rm (Proof)  From Theorem~\ref{th:w-choice}(I) we have  
${\bf C}_u(X\cup Y)\cap{\bf C}_u(C(X)\cup (Y\setminus X))\neq\emptyset$. 
Then under the present assumption we have 
$|{\bf C}_u(X\cup Y)\cap{\bf C}_u(C(X)\cup (Y\setminus X))|=1$,
which implies $C(X\cup Y)=C(C(X)\cup (Y\setminus X))$.
\pend}
\end{theorem}

Note that $C(X\cup Y)=C(C(X)\cup (Y\setminus X))$ for all $X, Y\in 2^E$
if and only if $C(X\cup Y)=C(C(X)\cup Y)$ for all $X, Y\in 2^E$. 
(This fact seems to be a folklore, but we give its proof 
in Appendix~\ref{Appendix2}.)\, 
Hence Theorem~\ref{th:choice2} is equivalent to Theorem~\ref{th:HKYY2022}.

\subsection{A choice function associated with an ordinally weak-concave 
function}\label{sec:Umax}

For any function $u: 2^E \rightarrow \mathbb{R}$ and $U\in 2^E$ define 
\begin{align*}
{\bf C}_u^{-1}(U)=\left\{X \in 2^E \mid U\in {\bf C}_u(X)\right\}.
\end{align*}
Here it should be noted that ${\bf C}_u^{-1}$ is not the inverse of 
the mapping (choice correspondence) 
${\bf C}_u: 2^E\to 2^{2^E}$ in a mathematical sense. 

\begin{lemma}\label{lem:lex0}
Let $u: 2^E \rightarrow \mathbb{R}$ be an ordinally $w$-concave 
function. 
For any $U \in 2^E$ with ${\bf C}_u^{-1}(U)\neq \emptyset$,\, 
if $X, Y \in {\bf C}_u^{-1}(U)$,\, then we have 
$X \cup Y \in {\bf C}_u^{-1}(U)$.
\medskip\\
{\rm (Proof) 
Consider $U \in 2^E$ with ${\bf C}_u^{-1}(U)\neq \emptyset$ and 
$X, Y \in {\bf C}_u^{-1}(U)$. We show $X \cup Y \in {\bf C}_u^{-1}(U)$. 
We can suppose that $X \neq Y$.

Now suppose to the contrary that $X \cup Y \notin {\bf C}_u^{-1}(U)$. 
Choose any $V \in {\bf C}_u(X \cup Y)$ in such a way that the following 
($*$) holds:
\begin{itemize}
\item[($*$)] $V$ attains the minimum of $|V \Delta U|$.
\end{itemize}
Since $U\notin {\bf C}_u(X \cup Y)$, we have $V \neq U$. 
Because of the ordinal w-concavity of $u$, there exist distinct 
$x \in$ $(U \backslash V) \cup\{\emptyset\}$ and 
$y \in(V \backslash U) \cup\{\emptyset\}$ such that
\begin{itemize}
\item[(i)] $u(U)<u(U-x+y)$,\, or
\item[(ii)] $u(V)<u(V-y+x)$,\, or
\item[(iii)] $u(U)=u(U-x+y)$\, and\, $u(V)=u(V-y+x)$.
\end{itemize}
If (iii) holds, then because of the definition of $V$ and since 
$V-y+x \subseteq X \cup Y$, we have $V-y+x \in {\bf C}_u(X \cup Y)$. 
This contradicts the assumption ($*$) for the choice of $V$. 
Also Condition (ii) contradicts $V \in {\bf C}_u(X \cup Y)$. Moreover, 
if (i) holds, then we must have $y \in Y \backslash X$ due to the definition 
of $U$ for $X$, which then contradicts the definition of $U$ for $Y$. 
This completes the proof of the present lemma.
\pend}
\end{lemma}

Lemma~\ref{lem:lex0} implies the following theorem.

\begin{theorem}\label{prop:1}
Let $u: 2^E \rightarrow \mathbb{R}$ be an ordinally $w$-concave function.
Then, for any $U \in 2^E$ with ${\bf C}_u^{-1}(U)\neq \emptyset$, there 
uniquely exists a set $U^{+} \in$ $2^E$ such that 
${\bf C}_u^{-1}(U)=\left[U, U^{+}\right]$ being an interval of $2^E$.
\end{theorem}

For an ordinally $w$-concave function $u: 2^E \rightarrow \mathbb{R}$ 
that satisfies the unique-maximizer condition {\bf (UM)},  
let 
$C: 2^E\to 2^E$ be a choice function associated with $u$.
Let us call every $U\in{\rm Im}(C)\equiv \{C(X)\mid X\in 2^E\}$ 
a {\it choice-set} and 
$U^+$ the {\it enclosure} of $U$. Also we call $X\in 2^E$ a {\it proper set} 
with respect to $u$ if for every $x\in E\setminus X$,\, $C(X+x)\neq C(X)$. 
The enclosure $U^+$ of a choice-set $U\in{\rm Im}(C)$ is the unique
maximal proper set $X$ that contains $U$ and satisfies $C(X)=U$.

\begin{remark}{\rm 
Alva and Do\u{g}an \cite{AlvaDogan2021} have shown that for any path-independent
choice function $C: 2^E\to 2^E$ and $U\in 2^E$ 
there uniquely exists a set $U^+ \in 2^E$
such that for every $X\in 2^E$, 
$C(X)=U$ if and only if $U\subseteq X\subseteq U^+$. 
This fact also follows from Theorem~\ref{prop:1} since 
any path-independent choice function $C$ is associated with an 
ordinally concave function $u$ that satisfies the unique-maximizer condition 
{\bf (UM)} (due to Yokote et~al.~\cite{YHKY2023}). 
\pend}
\end{remark}

\section{Discussions}\label{sec:discussions}

\subsection{Duality in ordinal concavity}\label{sec:duality}

Consider any ordinally concave function $u: 2^E\to\mathbb{R}$. 
In the definition of ordinal concavity, Definition~\ref{def:1},
the choice of $X$ and $X'$ is for an unordered pair, while the 
choice of $x$ and $x'$ and the associated conditions are  
given for an ordered pair $(x,x')$.
If we change the roles of $X$ and $X'$ we have an equivalent definition
of ordinal concavity as follows.

\begin{definition}[{\bf Ordinal Concavity*}] \label{def:3}
A function $u: 2^E\to \mathbb{R}$ satisfies  
{\rm ordinal concavity} if\, for every $X, X'\in 2^E$\, 
the following statement holds\,{\rm :} \smallskip\\
For every $x'\in X'\setminus X$ there exists 
$x\in (X\setminus X')\cup\{\emptyset\}$
such that 
\begin{itemize}
\item[{\rm (i)}] $u(X)<u(X+x'-x)$,\, or 
\item[{\rm (ii)}] $u(X')<u(X'-x'+x)$,\, or
\item[{\rm (iii)}] $u(X)=u(X+x'-x)$\, and\, $u(X')=u(X'-x'+x)$.\vspace{-0.1cm}
\end{itemize} 
\end{definition}

We discuss some implications of this fact. 
For any given function $u: 2^E\to \mathbb{R}$ let us define 
$u^\bullet: 2^E\to\mathbb{R}$ by
\begin{equation}\label{eq:dd1}
 u^\bullet(X)=u(E\setminus X)\qquad (\forall X\in 2^E).
\end{equation}
We call such $u^\bullet$ the {\it dual} of $u$. We may consider that 
$u^\bullet$ is defined on the dual Boolean lattice of $2^E$.\, 
Note that $(u^\bullet)^\bullet=u$.

\begin{lemma}\label{lem:dd1}
A function $u: 2^E\to \mathbb{R}$ is ordinally concave 
if and only if its dual $u^\bullet: 2^E\to\mathbb{R}$ is ordinally 
concave.
\medskip\\
{\rm (Proof) We can see that Definition~\ref{def:3} for $u$ gives 
exactly Definition~\ref{def:1} for $u^\bullet$ by considering 
the complements of $X$ and $X'$.
\pend}
\end{lemma}

Moreover, we also have the following lemma for ordinal w-concavity.

\begin{lemma}\label{lem:dd2}
A function $u: 2^E\to \mathbb{R}$ is ordinally w-concave 
if and only if its dual $u^\bullet: 2^E\to\mathbb{R}$ is ordinally 
w-concave.
\medskip\\
{\rm (Proof) The definition of ordinal w-concavity, Definition~\ref{def:2}, 
is self-dual, so that the present lemma holds.
\pend}
\end{lemma}

Hence, as a metatheory, if we have a valid statement for an ordinally 
(w-)concave function $u$, then the statement obtained by dualization by 
taking complements is also valid. 

For example, as a dual of Theorem~\ref{th:w-choice}(I) we have 
the following theorem.

\begin{theorem}\label{th:w-choice2}
Suppose that $u: 2^E\to\mathbb{R}$ satisfies ordinal concavity. 
Then, for every $X, Y\in 2^E$ and every $U\in{\bf C}_u(X\cup Y,E)$ 
there exists $Z\in{\bf C}_u(X,E)$ such that 
\[U\setminus(X\cup Y)\subseteq Z\setminus(X\cup Y).\]
{\rm (Proof) 
Let $U'=E\setminus U\in{\bf C}_{u\bullet}(E\setminus(X\cup Y))$.
Since $u^\bullet$ is ordinally concave due to Lemma~\ref{lem:dd1},
it follows from Theorem~\ref{th:w-choice}(I) that putting 
$X\gets E\setminus(X\cup Y)$ and $X\cup Y\gets E\setminus X$, 
there exists $Z'\in{\bf C}_{u^\bullet}(E\setminus X)$ such that 
$Z'\cap(E\setminus(X\cup Y))\subseteq U'$.
Note that putting $Z=E\setminus Z'$,\, we have
\[Z'\cap(E\setminus(X\cup Y))\subseteq U' \quad \Longleftrightarrow\quad
  U\setminus(X\cup Y)\subseteq Z\setminus(X\cup Y).\]
This completes the proof of the present theorem. \pend}
\end{theorem}

It may be worth considering another function associated with $u$ as follows.
\begin{equation}\label{eq:dual2}
  u^{\#}(X)=u(E)-u(E\setminus X)\qquad (\forall X\in 2^E).
\end{equation}
It should be noted that $u^{\#}: 2^E\to\mathbb{R}$ is 
{\sl ordinally {\rm (}w-{\rm )}convex} (i.e., $-u^{\#}$ is ordinally 
(w-)concave) when $u$ is 
ordinally (w-)concave. We have  $u^{\#}(X)=u(E)-u^\bullet(X)$ for all 
$X\in 2^E$ and we may call $u^{\#}$ the 
{\it dual ordinally {\rm (}w-{\rm )}convex function} 
of the ordinally (w-)concave function $u$.
Also note that when $u(\emptyset)=0$, we have $(u^{\#})^{\#}=u$.

\subsection{Domains of functions}\label{sec:M}

All the functions considered above have the unit hypercube or 
Boolean lattice $2^E$ as their domains.
Let us consider any sets as domains instead and examine 
how our above arguments work for the new problem setting. 

Let ${\bf Q}$ be a nonempty subset of $2^E$ 
and 
consider a function $u: {\bf Q}\to\mathbb{R}$. 
(Formally we may also
consider $u: 2^E\to\mathbb{R}\cup\{-\infty\}$ by putting $u(X)=-\infty$ for all
$X\in 2^E\setminus {\bf Q}$ so that ${\bf Q}$ is the effective domain 
${\rm dom}(u)=\{X\in 2^E\mid u(X)>-\infty\}$.)
\begin{itemize}
\item We call $u$ an {\it ordinally concave function} on ${\bf Q}$ if 
it satisfies Definition~\ref{def:1} with $2^E$ being replaced by ${\bf Q}$.
\item Also we call $u$ an {\it ordinally w-concave function} on ${\bf Q}$ if 
it satisfies Definition~\ref{def:2} with $2^E$ being replaced by ${\bf Q}$.
\end{itemize}
Here it should be noted that in the definition of ordinal (w-)concavity
we put $u(X)=-\infty$ for every $X\in 2^E\setminus {\bf Q}$ and we define
$-\infty < \alpha$\, for all\, $\alpha\in\mathbb{R}$. 
Then we see the following facts.
\begin{enumerate}
\item In Section~\ref{sec:w-concavity0} define a minor of $u$ as follows: 
for any $X, Y\in{\bf Q}$ such that $X\subset Y$ define 
${\bf Q}^Y_X=\{Z\in 2^{Y\setminus X}\mid X\cup Z\in{\bf Q}\}$ and 
$u^Y_X(Z)=u(Z\cup X)-u(X)$ for each $Z\in {\bf Q}^Y_X$.
\item All the statements given in Sections~\ref{sec:w-concavity}, 
\ref{sec:characterization}, and \ref{sec:strengthening} hold true 
by replacing $2^E$ by ${\bf Q}$ with $\emptyset\in {\bf Q}$.
\item 
In Section~\ref{sec:choice2} suppose $\emptyset \in \mathbf{Q}$ and define 
a choice function $C: 2^E \rightarrow \mathbf{Q}$ in such a way that 
$C(X) \in 2^X \cap \mathbf{Q}$ for each $X \in 2^E$. Also, define the choice 
correspondence rationalized by $u$ in such a way that 
$\mathbf{C}_u(X)=\operatorname{Arg} \max \left\{u(Z) \mid Z \in 2^X \cap \mathbf{Q}\right\}$ 
for each $X \in 2^E$. If $C(X) \in \mathbf{C}_u(X)$ for each $X \in 2^E$, 
we call $C$ a choice function associated with $u$. Then all the statements 
in Section~\ref{sec:choice2} hold true {\it mutatis mutandis.}
\item If $u: {\bf Q}\to\mathbb{R}$ is ordinally w-concave and is a constant
function on ${\bf Q}$, then ${\bf Q}$ is an M${}^\natural$-convex set.
\end{enumerate}

\subsection{A lexicographic composition of two functions}\label{sec:lex}

Let us consider the lexicographical order $\le_\ell$ on $\mathbb{R}^2$ 
defined by 
$(a,b)<_\ell (c,d)
\ \ \Longleftrightarrow\ \ 
{\rm (i)\ } a<c\ \ {\rm or\ \ (ii)\ }\ a=c\ {\rm \ and\ } b<d$,\, 
for all $a, b, c, d\in\mathbb{R}$.
Let $(\mathbb{R}^2)_\ell$ be the set $\mathbb{R}^2$ endowed with the 
lexicographical order $\le_\ell$. 

Consider two functions $u_i: 2^E\to\mathbb{R}$ for $i=1,2$.
For any $X\in 2^E$ let $(\hat{X}_1,\hat{X}_2)$ be the lexicographic
 maximizer $(X_1,X_2)$ of 
\begin{equation}\label{eq:1}
{\rm Lexico}\max\{(u_1(X_1),u_2(X_2))\in\mathbb{R}^2\mid 
   X_1,X_2\subseteq X,\, X_2=X\setminus X_1\}.
\end{equation}
Then define a function $\hat{u}: 2^E\to(\mathbb{R}^2)_\ell$ 
as follows. 
\begin{equation}\label{eq:2}
  \hat{u}(X)=(u_1(\hat{X}_1),u_2(\hat{X}_2))
\in(\mathbb{R}^2)_\ell
\qquad (\forall X\in 2^E). 
\end{equation}
We write $\hat{u}(X)=(\hat{u}_1(X),\hat{u}_2(X))$.
Note that for each $X\in 2^E$ we have 
\begin{equation}\label{eq:3}
\hat{X}_1\in {\rm Arg}\max\{u_2(X\setminus Z)\mid Z\in {\bf C}_{u_1}(X)\}
\end{equation}
and
\begin{equation}\label{eq:3a}
 \hat{u}(X)=(\hat{u}_1(X),\hat{u}_2(X))
 =(u_1(\hat{X}_1),u_2(X\setminus\hat{X}_1)).
\end{equation}
We call $\hat{u}$ the {\it lexicographic composition} of the ordered pair 
$(u_1,u_2)$ of 
functions $u_1$ and $u_2$. 
Let us denote $\hat{u}=u_1\diamond u_2$. 
The lexicographic composition $u_1\diamond u_2$
of $u_1$ and $u_2$ is ordinally w-concave if 
for every $X, Y\in 2^E$ with $X\neq Y$\, 
the following statement holds\,{\rm :} \smallskip\\
There exist distinct 
$x\in (X\setminus Y)\cup\{\emptyset\}$ and 
$y\in (Y\setminus X)\cup\{\emptyset\}$
such that 
\begin{itemize}
\item[{\rm (i)${}_\ell$}] $\hat{u}(X)<_\ell \hat{u}(X-x+y)$\,,\, {\rm or} 
\item[{\rm (ii)${}_\ell$}] $\hat{u}(Y)<_\ell \hat{u}(Y-y+x)$\,,\, {\rm or}
\item[{\rm (iii)${}_\ell$}] $\hat{u}(X)=\hat{u}(X-x+y)$\, {\rm and}\, 
$\hat{u}(Y)=\hat{u}(Y-y+x)$\,.\vspace{-0.1cm}
\end{itemize} 

The lexicographic composition $\hat{u}=u_1\diamond u_2$ can be regarded 
as an ordinal analogue of convolution in convex analysis. 
However, 
the lexicographic composition $\hat{u}=u_1\diamond u_2$ is not 
ordinally w-concave even if both $u_1$ and $u_2$ are ordinally w-concave, 
 in general. (See an example confirming this claim 
in Appendix~\ref{Appendix3}.)

We have the following theorem on the lexicographic composition 
$u_1\diamond u_2$ for a special class of ordinally concave 
functions $u_1$.

\begin{theorem}\label{th:lex}
If $u_1: 2^E\to\mathbb{R}$ is an ordinally w-concave function 
 that satisfies the unique-maximizer condition {\bf (UM)}
and $u_2: 2^E\to\mathbb{R}$ is ordinally w-concave, then 
 the lexicographic composition $u_1\diamond u_2$
is ordinally w-concave.
\medskip\\
{\rm (Proof) 
For any $X,Y\in 2^E$ with $X\neq Y$ 
consider $\hat{X}\in 2^X$ and $\hat{Y}\in 2^Y$ 
satisfying
\begin{equation}\label{eq:4*}
\hat{X}\in {\bf C}_{u_1}(X), \qquad \hat{Y}\in{\bf C}_{u_1}(Y).
\end{equation}
\noindent
{\bf Case 1}: Suppose that $\hat{X}\neq\hat{Y}$.
By ordinal w-concavity of $u_1$ there exist distinct 
$\hat{x}\in (\hat{X}\setminus \hat{Y})\cup\{\emptyset\}$ and 
$\hat{y}\in (\hat{Y}\setminus \hat{X})\cup\{\emptyset\}$
such that 
\begin{itemize}
\item[{\rm (i)${}_1$}] ${u_1}(\hat{X})< {u_1}(\hat{X}-\hat{x}+\hat{y})$\,,\, or 
\item[{\rm (ii)${}_1$}] ${u_1}(\hat{Y})< {u_1}(\hat{Y}-\hat{y}+\hat{x})$\,,\, or
\item[{\rm (iii)${}_1$}] ${u_1}(\hat{X})={u_1}(\hat{X}-\hat{x}+\hat{y})$\, 
and\, 
${u_1}(\hat{Y})={u_1}(\hat{Y}-\hat{y}+\hat{x})$\vspace{-0.1cm}
\end{itemize}
Suppose that (i)${}_1$ holds. 
Then, because of the definition of $\hat{X}$ we have 
$\hat{X}-\hat{x}+\hat{y}\not\subseteq X$, which implies 
$\hat{y}\in Y\setminus X$ and 
$\hat{X}-\hat{x}+\hat{y}\subseteq X+\hat{y}$. It follows from (i)${}_1$ that
\begin{equation}\label{eq:5}
\hat{u}_1(X)=u_1(\hat{X})< {u_1}(\hat{X}-\hat{x}+\hat{y})
\le\hat{u}_1(X+\hat{y}),
\end{equation}
where recall the notation of $\hat{u}_1$ (and $\hat{u}_2$) in (\ref{eq:3a}).
Hence we see that Condition~(i)${}_\ell$ holds 
for $x=\emptyset$ and $y=\hat{y}\in Y\setminus X$. 

Similarly, if (ii)${}_1$ holds, then we can show that 
Condition (ii)${}_\ell$ holds
for $x=\hat{x}\in X\setminus Y$ and $y=\emptyset$.

Suppose that Condition (iii)${}_1$ holds. Since $u_1$ satisfies the 
unique-maximizer condition ${\bf (UM)}$, we have 
$\hat{X}-\hat{x}+\hat{y}\not\subseteq X$, which implies 
$\hat{y}\in Y\setminus X$ and $\hat{X}-\hat{x}+\hat{y}\subseteq X+\hat{y}$. 
Since $\hat{X}$ and $\hat{X}-\hat{x}+\hat{y}$ are distinct subsets of 
$X+\hat{y}$ with ${u_1}(\hat{X})={u_1}(\hat{X}-\hat{x}+\hat{y})$, 
by ${\bf (UM)}$ again there uniquely exists 
$\hat{X}'\in {\bf C}_{u_1}(X+\hat{y})$ such that 
\begin{equation}
 u_1(\hat{X}')> u_1(\hat{X})=u_1(\hat{X}-\hat{x}+\hat{y}).
\end{equation}
Therefore,
\begin{equation}
 \hat{u}_1(X)=u_1(\hat{X})<u_1(\hat{X}')=\hat{u}_1(X+\hat{y}).
\end{equation}
Hence we see that Condition (i)${}_\ell$ holds for $x=\emptyset$ and 
$y=\hat{y}\in Y\setminus X$.
\smallskip\\
{\bf Case 2}: Suppose that $\hat{X}=\hat{Y}$. Put $U=\hat{X}(=\hat{Y})$,
 $\tilde{X}=X\setminus U$ and $\tilde{Y}=Y\setminus U$. 
Since $\tilde{X}\neq \tilde{Y}$ and $u_2$ is 
ordinally w-concave, there exist distinct 
$\tilde{x}\in (\tilde{X}\setminus\tilde{Y})\cup\{\emptyset\}$ and 
$\tilde{y}\in (\tilde{Y}\setminus\tilde{X})\cup\{\emptyset\}$
such that
\begin{itemize}
\item[{\rm (i)${}_2$}] 
$u_2(\tilde{X})
< u_2(\tilde{X}-\tilde{x}+\tilde{y})$\,,\, or 
\item[{\rm (ii)${}_2$}] 
$u_2(\tilde{Y})
< u_2(\tilde{Y}-\tilde{y}+\tilde{x})$\,,\, or 
\item[{\rm (iii)${}_2$}] 
$u_2(\tilde{X})
=u_2(\tilde{X}-\tilde{x}+\tilde{y})$\, and\, 
$u_2(\tilde{Y})
=u_2(\tilde{Y}-\tilde{y}+\tilde{x})$\,.\vspace{-0.1cm}
\end{itemize} 

If $\tilde{x}\neq\emptyset$, we have $\tilde{x}\in \tilde{X}\subseteq X$. 
Since $\tilde{x}\notin\tilde{Y}=Y\setminus U$, we have $\tilde{x}\notin Y$.
Similarly we can show that $\tilde{y}\notin X$ if $\tilde{y}\neq\emptyset$. 
Consequently, we have $\tilde{x}\in (X\setminus Y)\cup\{\emptyset\}$ 
and $\tilde{y}\in(Y\setminus X)\cup\{\emptyset\}$.
Also note that ${X}-\tilde{x}+\tilde{y}, {Y}-\tilde{y}+\tilde{x}
\supseteq U=\hat{X}=\hat{Y}$.\smallskip

\noindent
{\bf Case 2(i)}: Suppose that (i)${}_2$ holds. 
Then, since $U\subseteq {X}-\tilde{x}+\tilde{y}\subseteq X\cup Y$, 
it follows from Lemma~\ref{lem:lex0} and (i)${}_2$ that we have 
\begin{equation}\label{eq:x2}
\hat{u}_1({X})=\hat{u}_1({X}-\tilde{x}+\tilde{y}),\quad
\hat{u}_2({X})={u}_2(\tilde{X})<{u}_2(\tilde{X}-\tilde{x}+\tilde{y})
= \hat{u}_2({X}-\tilde{x}+\tilde{y}).
\end{equation}
Hence, for $x=\tilde{x}$ and 
$y=\tilde{y}$\, Condition (i)${}_\ell$ holds.\smallskip

\noindent
{\bf Case 2(ii)}:  When (ii)${}_2$ holds, we can similarly show that 
(ii)${}_\ell$ holds. \smallskip

\noindent
{\bf Case 2(iii)}: Suppose that (iii)${}_2$ holds. 
Then we have
\begin{equation}\label{eq:x3}
\hat{u}_1({X})=\hat{u}_1({X}-\tilde{x}+\tilde{y}),\quad
\hat{u}_2({X})={u}_2(\tilde{X})={u}_2(\tilde{X}-\tilde{x}+\tilde{y})
= \hat{u}_2({X}-\tilde{x}+\tilde{y}),
\end{equation}
\begin{equation}\label{eq:x4}
\hat{u}_1({Y})=\hat{u}_1({Y}-\tilde{y}+\tilde{x}),\quad
\hat{u}_2({Y})={u}_2(\tilde{Y})={u}_2(\tilde{Y}-\tilde{y}+\tilde{x})
= \hat{u}_2({X}-\tilde{y}+\tilde{x}).
\end{equation}
Hence (iii)${}_\ell$ holds.\smallskip

This completes the proof of the present theorem.
\pend}
\end{theorem}

\begin{remark}{\rm
If $u_1: 2^E\to\mathbb{R}$ is ordinally (w-)concave and 
$u_2: 2^E\to\mathbb{R}$ is M${}^\natural$-concave, then in order to compute
$\hat{u}(X)=u_1\diamond u_2(X)$ for $X\in 2^E$\, the lexicographic 
maximization of (\ref{eq:1}) is reduced to 
\begin{enumerate}
\item Maximization of $u_1^X$ to obtain ${\bf C}_{u_1}(X)$ and
\item Maximization of $(u_2^X)^\bullet: {\bf C}_{u_1}(X)\to \mathbb{R}$.
\end{enumerate}
The latter is a special case of maximization of the sum of two 
M${}^\natural$-concave functions, which can be solved 
if ${\bf C}_{u_1}(X)$ is 
appropriately identified for the maximization.
In particular, if $u_2$ is a modular function, i.e., 
$u_2(X)=\sum_{x\in X}u_2(x)$ for all $X\in 2^E$, then the latter maximization
can be solved by a greedy algorithm over ${\bf C}_{u_1}(X)$. 
\pend}
\end{remark}

\section{Concluding Remarks}\label{sec:concluding}

We have investigated combinatorial structures of ordinally concave functions 
and (newly introduced) ordinally w-concave functions. 
We have revealed their fundamental properties and facts such as \vspace{-0.1cm}
\begin{enumerate}
\item The local optimality implies the global optimality.\vspace{-0.3cm}
\item The set of maximizers in any interval of $2^E$ forms an 
M${}^\natural$-convex set.\vspace{-0.3cm}
\item We have shown that the above two properties characterize 
ordinal w-concavity.
\vspace{-0.3cm}
\item We have given an ${\rm O}(|E|^2)$ algorithm for maximizing 
ordinally concave functions with a function evaluation oracle.
\vspace{-0.3cm}
\item We have shown the duality in ordinal (w-)concavity and its implications.
\vspace{-0.3cm}
\item We have proposed the lexicographic composition of two ordinally 
w-concave functions.
\vspace{-0.1cm}
\end{enumerate}
It is worth further investigating the structures of ordinally (w-)concave 
functions in view of economics and discrete optimization. 
An algorithmic open problem is to maximize an ordinally w-concave function
in polynomial time, using a function evaluation oracle, even for some 
special class of ordinally w-concave functions.
Also, besides Theorem~\ref{th:lex} it is interesting to investigate any other
appropriate conditions for the lexicographic composition $u_1\diamond u_2$ 
of functions $u_1$ and $u_2$ to become an ordinally w-concave function.

\section*{Acknowledgments}

We are grateful to Kazuo Murota for his helpful comments on an earlier 
version of the present paper.
This work was supported by JST ERATO Grant Number JPMJER2301, Japan.
S.~Fuji\-shige's research was supported by JSPS KAKENHI Grant Numbers 
JP19K11839 
and JP22K11922 and by the Research Institute for Mathematical Sciences, 
an International Joint Usage/Research Center located in Kyoto University. 
F.~Kojima's research was supported by JSPS KAKENHI Grant Number JP21H04979.

\appendix

\section{Appendix}

\subsection{An example related to Theorem~\ref{th:w-choice}}
\label{Appendix:example3}

Let $E=\{a, b, c\}$ and define a function $u: 2^E \rightarrow \mathbb{R}$ 
as follows.
\begin{equation}\nonumber
\begin{array}{rccccccccc}
X= & \emptyset & \{a\} & \{b\} & \{c\} & \{a, b\} & \{a, c\} & \{b, c\} & \{a, b, c\} \\
\hline u(X)= & 0 & 4 & 3 & 1 & 2 & 5 & 6 & 0 
\end{array}
\end{equation}
We can see that the instance of $X=\{a\}, X^{\prime}=\{b,c\}$, and 
$x=a \in X \backslash X^{\prime}$ violates the definition of ordinal concavity. 
On the other hand, $u$ satisfies the definition of ordinal w-concavity.

Take $X=\{a,b\}$,\, $Y=\{c\}$,\, and\, $\{a\}\in \mathbf{C}_u(X)=\{\{a\}\}$,\, 
$\{b,c\}\in \mathbf{C}_u(X\cup Y)=\{\{b,c\}\}$.\, It holds that
\begin{align*}
\{b,c\}\cap \{a,b\}=\{b\}\nsubseteq \{a\}.
\end{align*}
Therefore, neither (I) nor (II) of Theorem~\ref{th:w-choice} holds.

\subsection{A proof}\label{Appendix2}

\begin{proposition} 
Let $u: 2^E\to\mathbb{R}$ be any ordinally concave function and $C: 2^E\to 2^E$
be a choice function associated with $u$. 
The following statements {\rm (a)} and {\rm (b)} are equivalent\,{\rm :}
\begin{itemize}
\item[{\rm (a)}] $C(X\cup Y)=C(C(X)\cup(Y\setminus X))$\, for all $X, Y\in 2^E$.
\item[{\rm (b)}] $C(X\cup Y)=C(C(X)\cup Y)$\, for all $X, Y\in 2^E$.
\end{itemize}
{\rm 
(Proof) The implication (b) $\Rightarrow$ (a) clearly holds. 
We show the converse.
Suppose that (a) holds. Consider any $X,Y\in 2^E$ and define 
\[Z=X\setminus((X\setminus C(X))\cap Y).\]
Then it follows from (a) that we have
\[C(X\cup Y)=C(Z\cup Y)=C(C(Z)\cup(Y\setminus C(X)))
=C(C(X)\cup(Y\setminus C(X)))=C(C(X)\cup Y),\]
where the second equality follows from (a) and $Y\setminus Z=Y\setminus C(X)$, 
the third equality follows from the fact that $C$ is associated with 
$u$, because $C(X)\subseteq Z\subseteq X$ and hence 
$C(Z)=C(X)$.
\pend}
\end{proposition}

\subsection{A lexicograhic composition of ordinally w-concave functions 
is not ordinally w-concave in general}\label{Appendix3}

Let $E=\{a,b,c,d\}$ and define functions $u_i: 2^E\to\mathbb{R}$ $(i=1,2)$ 
as follows.
{\small
\begin{equation}\nonumber
\begin{array}{rccccccccccc}
X= & \emptyset & \{a\} & \{b\} & \{c\} & \{d\}&\{a,b\}&\{a,c\}& \{a,d\}&
\{b,c\}&\{b,d\}&\{c,d\}\\
\hline
u_1(X)= & 0 & 2 & 2 & 8 & 8 & 2 & 3 & 4 & 5&6&7\\
\hline
u_2(X)= & 0 & 6 & 6 & 5 & 1 & 7 & 1 & 2 & 3&4&0\\
\end{array}
\end{equation}
\begin{equation}\nonumber
\begin{array}{rccccc}
X= & \{a,b,c\}
&\{a,b,d\}&\{a,c,d\}&\{b,c,d\}&\{a,b,c,d\}\\
\hline
u_1(X)= &1&1&1&1&0\\
\hline
u_2(X)= &1&5&0&0&0\\
\end{array}
\end{equation}
}
We can see that both $u_1$ and $u_2$ are ordinally w-concave. 
Consider $X=\{a,b,c,d\}$ and $Y=\{c,d\}$, where note that  
$(X\setminus Y)\cup\{\emptyset\}=\{a,b\}\cup\{\emptyset\}$ and 
$(Y\setminus X)\cup\{\emptyset\}=\{\emptyset\}$. 
Then, for the present $X$ and $Y$ we have that 
for any $x\in \{a,b\}$ and $y=\emptyset$ 
none of (i)${}_\ell$, (ii)${}_\ell$, and (iii)${}_\ell$ hold. 
Hence the lexicographic composition $u_1\diamond u_2$ is not 
ordinally w-concave.


\begin{thebibliography}{99}

\bibitem{AlvaDogan2021} S. Alva and B. Do\u{g}an: Choice and market design.
In: {\it Online and Matching-Based Market Design} 
(F. Echenique, N. Immorlica, and V. V. Vazirani, Eds.; 
Cambridge University Press, 2023), pp.~238--263.
arXiv:2110.15446v2 [econ.TH] 4 Nov 2021.


\bibitem{CheKimKojima2021} Y.-K.~Che, J. Kim, and F. Kojima: 
Weak monotone comparative statistics. 
arXiv:1911.06442v4 [econ.TH] 24 Nov 2021.
\bibitem{ChenLi2021} X. Chen and M. Li: M${}^\natural$-convexity and its 
applications in operations. {\it Operations Research} {\bf 69} (2021) 
1396--1408.

\bibitem{FarooqShioura2005} R. Farooq and A. Shioura: A note on the 
equivalence between substitutability and M${}^\natural$-convexity. 
{\it Pacific Journal of Optimization} {\bf 1} (2005) 243--252. 

\bibitem{Fleiner2003} T. Fleiner: A fixed-point approach to stable matchings 
and some applications. {\it Mathematics of Operations Research} {\bf 28} 
(2003) 103--126.

\bibitem{Frank84} A. Frank: Generalized polymatroids. 
In: {\it Finite and Infinite Sets, I}\, (A.~Hajnal, L.~Lov\'{a}sz and 
V.~T.~S\'{o}s, 
eds., Colloquia Mathematica Socientatis J\'{a}nos Bolyai {\bf 37}, 
North-Holland, 1984), pp.~285--294.

\bibitem{Fuji2015} J. Fuji: Substitutable choice functions and convex 
geometry. {\it Discrete Applied Mathematics} {\bf 180} (2015) 283--285.

\bibitem{Fuji2005} S. Fujishige: {\it Submodular Functions and Optimization}
(Second Edition) (Annals of Discrete Mathematics {\bf 58}, Elsevier, 2005). 

\bibitem{HKYY2022}  I.~E.~Hafalir, F.~Kojima,  M.~B.~Yenmez, and K.~Yokote:
Design on matroids: Diversity vs.~meritocracy.
arXiv:2301.00237v1 [econ.TH] 31 Dec 2022.

\bibitem{Hassin82} R. Hassin: Minimum cost flow with set-constraints. 
{\it Networks} {\bf 12} (1982) 1--21.

\bibitem{HatfieldMilgrom2005} J. W. Hatfield and P. Milgrom: 
Matching with contracts. {\it American Economic Review} {\bf 95} (2005) 
913--935.

\bibitem{Koshevoy1999} G. A. Koshevoy: Choice functions and abstract 
convex geometries. {\it Mathematical Social Sciences} {\bf 38} (1999) 35--44.

\bibitem{MR2001} B. Monjardet and V. Raderanirina: The duality between 
the anti-exchange closure operators and the path independent choice 
operators on a finite set. {\it Mathematical Social Sciences} {\bf 41} 
(2001) 131--150.

\bibitem{Murota2003} K. Murota: {\it Discrete Convex Analysis} 
(SIAM Monographs on Discrete Mathematics and Applications {\bf 10}, 
SIAM, 2003).

\bibitem{Murota2016} K. Murota: Discrete convex analysis: A tool for 
economics and game theory. {\it Journal of Mechanism and Institution Design}
 {\bf 1} (2016) 151--273. Also see: Revised version arXiv:2212.03598v1 
[math.CO] 7 Dec 2022.

\bibitem{MurotaShioura1999} K. Murota and A. Shioura: M-convex function
on generalized polymatroid. {\it Mathematics of Operations Research} 
{\bf 24} (1999) 95--105.

\bibitem{MurotaShioura2003} K. Murota and A. Shioura: Quasi M-convex and 
L-convex functions: quasi-convexity in discrete optimization. 
{\it Discrete Applied Mathematics} {\bf 131} (2003) 467--494.

\bibitem{MurotaShioura2024} K. Murota and A. Shioura: Note on minimization 
of quasi M${}^\natural$-convex functions. 
{\it Japan Journal of Industrial and Applied Mathematics} {\bf 41} (2024) 
857--880.

\bibitem{ShiouraTamura2015} A. Shioura and A. Tamura: Gross substitutes 
condition and discrete concavity for multi-unit valuations: a survey.
{\it Journal of the Operations Research Society of Japan} {\bf 58} 
(2015) 61--103.

\bibitem{Sotomayor1999} M. Sotomayor: Three remarks on the many-to-many stable 
matching problem. {\it Mathematical Social Sciences} {\bf 38} (1999) 55--70.

\bibitem{YHKY2023} K.~Yokote, I.~E.~Hafalir, F.~Kojima, and M.~B.~Yenmez: 
Rationalizing path-independent choice rules.\,
arXiv:2303.00892v2 [econ.TH] 28 May 2024.


\end{thebibliography}
\end{document}